\documentclass[12pt]{article}

\title{Retrieving information from subordination }

\author{Jean Bertoin\thanks{Laboratoire de Probabilit\'es et Mod\`eles Al\'eatoires, 
UPMC, 175 rue du Chevaleret, 75013 Paris, France. Email: jean.bertoin@upmc.fr} 
\and
Marc Yor \thanks{Laboratoire de Probabilit\'es et Mod\`eles Al\'eatoires, 
UPMC, 175 rue du Chevaleret, 75013 Paris, France. Email: deaproba@proba.jussieu.fr} 
\thanks{Institut Universitaire de France.}
 }

\date{ }

\usepackage{amsfonts,amsmath,amssymb,bbm}
\usepackage{setspace} 
\def\QED{\hfill $\Box$}

\font\tenmath=msbm10 scaled 1200
\font\sevenmath=msbm7 scaled 1200
\font\Fivemath=msbm5 scaled 1200
\newfam\mathfam \textfont\mathfam=\tenmath
\scriptfont\mathfam=\sevenmath \scriptscriptfont\mathfam=\Fivemath

\def \\ { \cr }
\def\R{\mathbb{R}}
\def \1{1 \mkern -6mu 1} 
\def\N{\mathbb{N}}
\def\E{\mathbb{E}}
\def\P{\mathbb{P}}
 
\def\f{\mathcal{F}}
\def\Q{\mathbb{Q}}
\def\R{\mathbb{R}}

\def \d{{\rm d}}
\def \Y{{\mathcal Y}}
\def \X{{\mathcal X}}
\def \H{{\mathcal H}}
\def \I{{\mathcal I}}
\def \a{{\mathcal A}}

\newtheorem{theorem}{Theorem}

\newtheorem{proposition}{Proposition}

\begin{document}

\maketitle

\begin{abstract}  We show that if  $(X_s, s\geq 0)$ is a right-continuous process, $Y_t=\int_0^t\d s X_s$ its integral process and
$\tau = (\tau_{\ell}, \ell \geq 0)$ a subordinator, then the time-changed process  $(Y_{\tau_{\ell}}, \ell\geq 0)$ allows to retrieve the information
about $(X_{\tau_{\ell}}, \ell\geq 0)$ when  $\tau$ is stable, but not when $\tau$ is a gamma subordinator.
This question has been motivated by a striking identity in law involving the Bessel clock taken at an independent inverse Gaussian variable.
\end{abstract}


\begin{section}{Introduction and main statements}

\subsection{Motivation}
In Dufresne and Yor \cite{DY}, it was remarked that by combining Bougerol's identity in law (see, e.g. Bougerol \cite{Bougerol} and Alili {\it et al.} \cite{ADY}) and the symmetry principle of D\'esir\'e Andr\'e, there is the identity in distribution for every fixed $\ell\geq 0$
\begin{equation}\label{E0}
H_{\tau_\ell} \stackrel{\mbox{\tiny (law)}}{=} \tau_{a(\ell)},
\end{equation}
where $a(\ell)= {\mathrm{Argsinh}}(\ell)=\log\left(\ell+\sqrt{1+\ell^2}\right)$, 
$$H_t=\int_0^t \d s R^{-2}_s\,,\qquad t\geq 0\,,$$
is the Bessel clock constructed from a two-dimensional Bessel process $(R_s, s\geq 0)$ issued from $1$, and $(\tau_{\ell}, \ell\geq 0)$ is a stable $(1/2)$ subordinator independent from $(R_s, s\geq 0)$. 

In \cite{DY}, the authors wondered whether \eqref{E0} extends at the level of processes indexed by $\ell\geq 0$, or equivalently whether $(H_{\tau_\ell}, \ell \geq 0)$ has independent increments. Our main result entails that this is not the case. Indeed, Theorem \ref{T1} below  implies that for every $\ell\geq 0$, the filtration $(\hat \H_{\ell}, \ell\geq 0)$ generated by $(H_{\tau_\ell}, \ell\geq 0)$ contains the filtration generated by $\left(R_{\tau_{\ell}},  \ell\geq 0\right)$.
On the other hand,  $((R_s,H_s), s\geq 0)$ is a Markov (additive) process, and since subordination by an independent stable subordinator
preserves the Markov property,  $((R_{\tau_\ell},H_{\tau_\ell}), \ell \geq 0)$ is Markovian in its own filtration, which coincides with $\left(\hat \H_{\ell}, \ell \geq 0\right)$ by Theorem \ref{T1}. It is immediately seen that for any $\ell'>0$, the conditional distribution of $H_{\tau_{\ell+\ell'}}$ given $(R_{\tau_\ell},H_{\tau_\ell})$
does not only depend on $H_{\tau_\ell}$, but on  $R_{\tau_\ell}$ as well. Consequently the process $(H_{\tau_\ell}, \ell \geq 0)$ is not even an autonomous Markov process.
We point out that the process  $(R_{\tau_\ell}, \ell \geq 0)$ is Markov (by subordination), and refer to 
a forthcoming paper   \cite{BDY}  for details on the semigroup of $(R_{\tau_\ell},H_{\tau_\ell})$.

\subsection{Main results}
More generally, we consider in this work  an $\R^d$-valued process  $(X_s,s\geq 0)$ with right-continuous sample paths, and $(\tau_{\ell}, \ell\geq 0)$ a stable subordinator with index $\alpha\in(0,1)$. We stress that we do  {\it not require $X$ and  $\tau$ to be independent}. 
Introduce
$$Y_u=\int_0^u\d s X_s\,,\qquad u\geq 0\,,$$
and the right-continuous time-changed processes
 $$\hat X_\ell=X_{\tau_{\ell}}  \hbox{ and } \hat Y_{\ell}=Y_{\tau_{\ell}} ,\qquad \ell \geq 0\,.$$
We are interested in comparing the information embedded in the processes $\hat X$ and $\hat Y$, respectively. We state our main result. 
 
 \begin{theorem}\label{T1}  The right-continuous filtration $(\hat \Y_{\ell}, \ell\geq 0)$ generated by the process $\left( \hat Y_{\ell}, \ell \geq 0\right)$
  contains  the right-continuous filtration $(\hat \X_{\ell}, \ell\geq 0)$ generated by  $\left(\hat X_{\ell}, \ell \geq 0\right)$. 
 \end{theorem}

A perusal of the proof (given below in Section 2) shows that Theorem \ref{T1} can be extended to the case when it is only assumed that $\tau$ is a subordinator such that the tail of its L\'evy measure is regularly varying at $0$ with index $-\alpha$, which suggests that this result might hold for more general subordinators. On the other hand, if $(N_{\ell}, \ell\geq 0)$ is
any increasing step-process issued from $0$, such as for instance a Poisson process,  then the time-changed process $(Y_{N_{\ell}}, \ell\geq 0)$ stays at $0$ until the first jump time of $N$ which is strictly positive a.s. This readily implies that the germ-$\sigma$-field
$$\bigcap_{\ell>0}\sigma(Y_{N_{\lambda}}, \lambda \leq \ell)$$
is trivial, in the sense that every event of this field has probability either $0$ or $1$.
Focussing on subordinators with infinite activity, it is interesting to point
out that Theorem \ref{T1}  fails when one replaces the stable subordinator $\tau$ by  a gamma subordinator, as it can be seen from the following observation.

\begin{proposition} \label{P1}  Let $\gamma=(\gamma_t, t\geq 0)$ be a gamma-subordinator and $\xi$ a random variable with values in
$(0,\infty)$ which is independent of $\gamma$. Then the germ-$\sigma$-field
$$\bigcap_{t>0}\sigma(\xi\gamma_s, s\leq t)$$
is trivial.
\end{proposition}
We point out that Proposition \ref{P1} holds more generally when $\gamma$ is replaced by a subordinator with logarithmic singularity, also called of class $({\mathcal L})$, in the sense that the drift coefficient is zero and the L\'evy measure is absolutely continuous with density $g$ such that
$g(x) = g_0 x^{-1} + G(x)$ where $g_0$ is some strictly positive constant and $G: (0,\infty)\to \R$ a measurable function such that
$$\int_0^1|G(x)|\d x<\infty\ ,\ g(x)\geq 0\ ,\ \hbox{and} \ \int_1^{\infty} g(x)\d x<\infty\,.$$ 
Indeed, it has been shown by von Renesse {\it et al.} \cite{RYZ} that such subordinators enjoy a quasi-invariance property analogous to that of the gamma subordinator, and this is the key to Proposition \ref{P1}. 

It is natural to investigate a similar question in the framework of stochastic integration. For the sake of simplicity, we shall focus on the one-dimensional case.
We thus consider a real valued Brownian motion $(B_t, t\geq 0)$ in some filtration$(\f_t, t\geq 0)$ and an $(\f_t)$-adapted continuous process
$(X_t, t\geq 0)$, and consider the stochastic integral
$$I_t=\int_0^t X_s \d B_s\,,\qquad t\geq 0\,.$$
We claim the following.

\begin{proposition}\label{P2} Fix $\eta>0$ and assume that the sample paths of $(X_t, t\geq 0)$ are  H\"older-continuous with exponent $\eta$ a.s.
Suppose also that $(\tau_{\ell}, \ell\geq 0)$ is a stable subordinator of index $\alpha\in(0,1)$, which is independent of $\f_{\infty}$.
Then the right-continuous filtration $(\hat \I_{\ell}, \ell\geq 0)$ generated by the subordinate stochastic integral  $\left( \hat I_{\ell}=I_{\tau_{\ell}}, \ell \geq 0\right)$
  contains  the right-continuous filtration generated by  $\left(|X_{\tau_{\ell}}|, \ell \geq 0\right)$. 
\end{proposition}

The proofs of these statements are given in the next section.

\end{section}

\begin{section} {Proofs}

\subsection{Proof of Theorem 1.}
 We first observe that the proof can be reduced to showing that the germ-$\sigma$-field $\hat \Y_0$ contains the $\sigma$-field 
generated by $\hat X_0=X_0$. Indeed, let us take this for granted, fix $\ell>0$ and define
$\tau'_u=\tau_{\ell+u}-\tau_{\ell}$ and $X'_v=X_{v+\tau_\ell}$. Then $\tau'$ is again a stable$(\alpha)$ subordinator and
$X'$ a right-continuous process, and 
$$\hat Y_{\ell + u}-\hat Y_{\ell}= \int_0^{\tau'_u} \d v X'_v\,.$$
Hence $X'_0=\hat X_\ell$ is measurable with respect to $\hat \Y_\ell$, and our claim follows.

Thus we only need to verify that $X_0$ is $\hat \Y_0$-measurable.  In this direction, we shall use the following version
of the Law of Large Numbers for the jumps $\Delta \tau_s=\tau_s-\tau_{s-}$ of a stable subordinator.
Fix any $m>2/\alpha$ and introduce for any given $b\in \R$ and $\varepsilon >0$
$$N_{\varepsilon,b}={\rm Card}\{s\leq \varepsilon: b\Delta\tau_s> \varepsilon^m\}\,.$$
Note  that  $N_{\varepsilon,b}\equiv 0$ for $b\leq 0$.
For the sake of simplicity, we henceforth suppose that the tail of the L\'evy measure of $\tau$ is $x\mapsto x^{-\alpha}$, which induces no loss of generality. So for $b>0$, $N_{\varepsilon,b}$ is a Poisson variable with parameter
$$\varepsilon (\varepsilon^m/b)^{-\alpha}=b^{\alpha}\varepsilon^{1-m\alpha}\,.$$
Combining a standard argument based on the Borel-Cantelli lemma and Chebychev's inequality   with monotonicity, we get that for $\varepsilon=1/n$
\begin{equation}\label{E1}
\lim_{n\to\infty} n^{1-\alpha m}N_{1/n,b}= b^{\alpha}\qquad \hbox{for all $b>0$, almost-surely.}
\end{equation}

Let us assume that  the process $X$ is real-valued as the case of higher dimensions will then follow by considering coordinates.
Set
$$J_{\varepsilon}={\rm Card}\{s\leq \varepsilon : \Delta \hat Y_s>\varepsilon^m\}\,,$$
where as usual $\Delta \hat Y_s=\hat Y_s-\hat Y_{s-}$.
We note that
$$\Delta \hat Y_s-X_0\Delta \tau_s=\int_{\tau_{s-}}^{\tau_s} \d u (X_u-X_0)\,.$$
Hence if we set  $a_{\varepsilon}=\sup_{0\leq u \leq \tau_{\varepsilon}}|X_u-X_0|$, then
$$\left(X_0-a_{\varepsilon}\right) \Delta \tau_s\leq \Delta \hat Y_s 
\leq \left(X_0+a_{\varepsilon} \right) \Delta \tau_s\,,$$
from which we deduce 
$$N_{\varepsilon, X_0-a_{\varepsilon}}\leq J_{\varepsilon}
\leq N_{\varepsilon, X_0+a_{\varepsilon}}\,.$$
Since $X$ has right-continuous sample paths a.s., we have $\lim_{\varepsilon\to 0} a_{\varepsilon}=0$ a.s., and 
taking $\varepsilon = 1/n$, we now deduce from \eqref{E1} that 
$$\lim_{n\to\infty} n^{1-m\alpha} J_{\varepsilon}= (X_0^+)^{\alpha}\,.$$
Hence $X_0^+$ is $\hat\Y_0$-measurable, and the same argument also shows that $X_0^-$ is $\hat\Y_0$-measurable.
This completes the proof of our claim. \QED

\subsection{Proof of Proposition 1.}
 Let $\Omega$ denote the space of c\`adl\`ag paths $\omega: [0,\infty)\to \R_+$ endowed with
the right-continuous filtration
$(\a_t, t\geq 0)$ generated by the canonical process $\omega_t=\omega(t)$, and write $\Q$ for the law on $\Omega$ of the  process
$(\xi \gamma_t, t\geq 0)$. 

It is well known that for every $x>0$ and $t>0$, the distribution of the process $(x\gamma_s, 0\leq s \leq t)$ is absolutely continuous with respect to that of the gamma process $(\gamma_s, 0\leq s \leq t)$ with density
$$x^{-t} \exp\left((1-1/x)\gamma_t\right)\,.$$
Because $\xi$ and $\gamma$ are independent, this implies that for any event $\Lambda\in\a_r$ with $r<t$
$$\Q\left(\Lambda\right)= \E\left( \xi^{-t} \exp\left((1-1/\xi)\gamma_t\right)\, {\bf 1}_{\{\gamma\in\Lambda\}}\right)\,.$$

Observe that
$$\lim_{t\to 0+} \xi^{-t} \exp\left((1-1/\xi)\gamma_t\right) = 1 \qquad \hbox{a.s.}$$
and the convergence also holds in $L^1(\P)$ by an application of Scheff\'e's lemma (alternatively, one may also invoke the convergence of backwards martingales). We deduce that for every 
$\Lambda\in \f_0$, we have
$$\Q(\Lambda) = \P(\gamma\in\Lambda)$$
and the right-hand-side must be $0$ or $1$ because the gamma process fulfills the Blumenthal's $0$-$1$ law. \QED

\subsection{Proof of Proposition 2.}  The guiding line is similar to that of the proof of Theorem \ref{T1}. In particular it suffices to verify that  $|X_0|$ is measurable with respect to the germ-$\sigma$-field $\hat \I_0$. 

Because $B$ and $\tau$ are independent,  the subordinate Brownian motion $(\hat B_{\ell}=B_{\tau_{\ell}}, \ell\geq 0)$ is a symmetric stable L\'evy process with index $2\alpha$. With no loss of generality, we may suppose that the tail of its L\'evy measure $\Pi$ is given by 
$\Pi(\R\backslash [-x,x])=x^{-2\alpha}$. As a consequence, for every $m>2/\alpha$ and  $\varepsilon>0$ and $b\in\R$, if define 
$$N_{\varepsilon, b}={\rm Card}\{s\leq \varepsilon: |b\Delta \hat B_s|^2>\varepsilon^m\}\,,$$ 
then $N_{\varepsilon,b}$ is a Poisson variable with parameter
$|b|^{2\alpha}\varepsilon^{1-m\alpha}$, and this readily yields 
\begin{equation}\label{E2}
\lim_{n\to\infty} n^{1-\alpha m}N_{1/n,b}= |b|^{2\alpha}\qquad \hbox{for all $b\in\R$, almost-surely.}
\end{equation}

Next set
$$J_{\varepsilon}={\rm Card}\{s\leq \varepsilon : |\Delta \hat I_s|^2>\varepsilon^m\}\,,$$
where as usual $\hat I_s=I_{\tau_s}$, and observe that
\begin{equation}\label{E5}
\Delta \hat I_s=X_0\Delta \hat B_s+ (X_{\tau_{s-}}-X_0)\Delta \hat B_s+\int_{\tau_{s-}}^{\tau_s}   (X_u-X_{\tau_{s-}})  \d B_u\,.
\end{equation}
Recall the assumption that the paths of $X$ are H\"older-continuous with exponent $\eta>0$, so the 
$(\f_t)$-stopping time
$$T=\inf\left\{u> 0: \sup_{0\leq v<u} (u-v)^{-\eta}|X_u-X_v|^2> 1\right\}$$
is strictly positive a.s. In particular, if we write $\Lambda_{\varepsilon}=\{\tau_{\varepsilon}<T\}$, 
then $\P(\Lambda_{\varepsilon})$ tends to $1$ as $\varepsilon\to 0+$.

We fix $a>0$, consider 
$$K_{\varepsilon,a}={\rm Card}\left\{s\leq \varepsilon : \left| \int_{\tau_{s-}}^{\tau_s} (X_u-X_{\tau_{s-}})\d B_u \right |^2>a \varepsilon^m\right\}\,,$$
and claim that
\begin{equation}\label{E3}
\lim_{\varepsilon\to 0}\varepsilon^{\alpha m -1}\E(K_{\varepsilon,a}, \Lambda_{\varepsilon})=0\,.
\end{equation}
If we take \eqref{E3} for granted, then we can complete the proof by an easy adaptation of the argument in Theorem \ref{T1}. 
Indeed, we can then find a strictly increasing sequence of integers $(n(k), k\in\N)$ such that with probability one, 
for all rational numbers $a>0$
\begin{equation}\label{E7}
\lim_{k\to\infty}n(k)^{1-\alpha m }K_{1/n(k),a}=0\,.
\end{equation}
We observe from \eqref{E5} that for any $a\in(0,1/2)$, if $|\Delta \hat I_s|^2>\varepsilon^m$, then necessarily either
$$|X_0\Delta \hat B_s|^2>(1-2a)^2\varepsilon^m\,,$$
 or 
 $$|(X_{\tau_{s-}}-X_0)\Delta \hat B_s|^2>a^2 \varepsilon^m\,,$$
  or
$$\left| \int_{\tau_{s-}}^{\tau_s}   (X_u-X_{\tau_{s-}})  \d B_u\right|^2>a^2 \varepsilon^m\,.$$
As 
$$\lim_{\varepsilon\to 0+}\sup_{0\leq s \leq \varepsilon}|X_{\tau_{s-}}- X_0|=0\,,$$ 
this easily entails, using \eqref{E2} and \eqref{E7}, that
\begin{eqnarray*}
\limsup_{k\to\infty}n(k)^{1-\alpha m } J_{1/n(k)} &\leq&  \lim_{k\to\infty}n(k)^{1-\alpha m } N_{1/n(k), (1-2a)^{-1}|X_0|} \\
 &=& (1-2a)^{-2\alpha}|X_0|^{2\alpha}\,,\qquad \hbox{a.s.}
 \end{eqnarray*}
where the identity in the second line stems from \eqref{E2}. A similar argument also gives
$$
\liminf_{k\to\infty}n(k)^{1-\alpha m } J_{1/n(k)} \geq  (1+2a)^{-2\alpha}|X_0|^{2\alpha}\,,\qquad \hbox{a.s.,}
$$
and as $a$ can be chosen arbitrarily close to $0$, we conclude that
$$
\lim_{k\to\infty}n(k)^{1-\alpha m } J_{1/n(k)} =|X_0|^{2\alpha}\,,\qquad \hbox{a.s.}
$$
Hence $|X_0|$ is $\hat \I_0$-measurable. 

Thus we need to establish \eqref{E3}. As $\tau$ is independent of $\f_{\infty}$, we have by an application of Markov's inequality that for every $s\leq \varepsilon$
\begin{eqnarray*}
& &\P\left(\left| \int_{\tau_{s-}}^{\tau_s} (X_u-X_{\tau_{s-}})\d B_u \right |^2>a \varepsilon^m, \Lambda_{\varepsilon}\mid   \tau\right)\\
& &\leq \frac{1}{a\varepsilon^m}\int_0^{\Delta \tau_s} \d v v^{\eta}
 \ \leq \ \frac{(\Delta \tau_s)^{1+\eta}}{a\varepsilon^m}\,.
\end{eqnarray*}
It follows that
\begin{eqnarray*}
\E(K_{\varepsilon,a}, \Lambda_{\varepsilon})
&\leq&\E\left(\sum_{s\leq \varepsilon}\left( \frac{(\Delta \tau_s)^{1+\eta}}{a\varepsilon^m}\wedge 1\right)
\right)\\
&=&\varepsilon c \int_{(0,\infty)}\d x x^{-1-{\alpha}}\left( \frac{x^{1+\eta}}{a\varepsilon^m}\wedge 1\right) 
= O(\varepsilon^{1-\alpha m/(1+\eta)})\,,
\end{eqnarray*}
where  for the second line we used the fact that the L\'evy measure of $\tau$ is $c x^{-1-\alpha}\d x$ for some unimportant constant $c>0$. 
This establishes \eqref{E3} and hence completes the proof of our claim. \QED

\end{section}

  \end{document}